\documentclass[twoside,reqno,11pt]{fcaa} 
 
\usepackage{graphicx}
\usepackage{epsfig}
\usepackage{amsthm}
\usepackage{amsmath}
\usepackage{latexsym}
\usepackage{amsfonts}
\usepackage{amssymb}

\newtheoremstyle{theorem}
  {15pt}          
  {15pt}  
  {\sl}  
  {\parindent}
  {\sc}  
  {. }   
  { }    
  {}     
\theoremstyle{theorem}

\newtheorem{theorem}{Theorem}[section]

\newtheoremstyle{defi}
  {15pt}          
  {15pt}  
  {\rm}  
  {\parindent}     
  {\sc}  
  {. }    
  { }    
  {}     
\theoremstyle{defi}

 \def\proofname{\indent\rm P\ r\ o\ o\ f.}
 \def\proofend{\hfill$\Box$}
 \def\qed{\hfill$\Box$} 

\usepackage{bm}
\usepackage{tikz}
\usepackage{pgfplots}

 \usepackage{hyperref} 

 \def\theequation{\arabic{section}.\arabic{equation}}

 \def\themycopyrightfootnote{\vspace*{3pt}
 \copyright \, 2014\,  Petr N. Vabishchevich
   \par  \noindent pp. 1--22
   \hfill  \vspace*{-36pt}
  }

 \title[NUMERICAL SOLUTION OF \dots]
{NUMERICAL SOLUTION \\ [4pt] OF NONSTATIONARY PROBLEMS \\ [4pt] FOR A SPACE-FRACTIONAL DIFFUSION EQUATION}
 \author[P. N. Vabishchevich]{Petr N. Vabishchevich$^{1,2}$}

\begin{document}

 \vbox to 2.5cm { \vfill }


 \bigskip \medskip

 \begin{abstract}
An unsteady problem is considered for a space-fractional diffusion equation in a bounded domain.
A first-order evolutionary equation containing a fractional power of an elliptic operator of second order
is studied for general boundary conditions of Robin type.
Finite element approximation in space is employed. 
To construct approximation in time, regularized two-level schemes are used.
The numerical implementation is based on solving the equation with the fractional power of the elliptic operator using an auxiliary Cauchy problem for a pseudo-parabolic equation.
The results of numerical experiments are presented for a model two-dimensional problem.

\medskip


 \smallskip


 \end{abstract}

 \maketitle

 \vspace*{-16pt}



 \section{Introduction}\label{sec:1}

\setcounter{section}{1}
\setcounter{equation}{0}\setcounter{theorem}{0}

Nowadays, non-local applied mathematical models based on the use of fractional derivatives in time and space 
are actively discussed \cite{baleanu2012fractional,eringen2002nonlocal,kilbas2006theory}. 
Many models, which are used in applied physics, biology, hydrology and finance,
involve both sub-diffusion (fractional in time) and supper-diffusion (fractional in space) operators. 
Supper-diffusion problems are treated as evolutionary problems with a fractional power of an elliptic operator.
For example, suppose that in a bounded domain $\Omega$ on the set of functions 
$u(\bm x) = 0, \ \bm x \in \partial \Omega$, 
there is defined the operator $\mathcal{A}$: $\mathcal{A}  u = - \triangle u, \ \bm x \in \Omega$. 
We seek the solution of the Cauchy problem for the equation with a fractional power of an elliptic operator:
\[
 \frac{d u}{d t} + \mathcal{A}^\alpha u = f(t),
 \quad 0 < t \leq T, 
\] 
\[
 u(0) = u_0,
\] 
for the given $f(\bm x, t)$, $u_0(\bm x), \ \bm x \in \Omega$  with $0 < \alpha < 1$ 
using the notation $f(t) = f(\cdot,t)$. 

To solve numerically evolutionary equations of first order, as a rule, two-level difference schemes are used 
for approximation in time.
Investigation of stability for such schemes in the corresponding finite-dimensional
(after discretization in space) spaces is based on the general theory of operator-difference schemes
\cite{Samarskii1989,SamarskiiMatusVabischevich2002}. In particular, the backward Euler scheme and Crank-Nicolson scheme are unconditionally stable for a non-negative operator. 
As for one-dimensional problems for the space-fractional diffusion equation,
an analysis of stability and convergence for this equation was conducted in \cite{jin2014error} using
finite element approximation in space. A similar study for the Crank-Nicolson scheme was considered earlier
in \cite{tadjeran2006second} using finite difference approximations in space.  

In discussing the problems of the numerical solution of multidimensional problems for the space-fractional diffusion equation,
emphasis is on spatial approximations. Many researchers (see, e.g., 
\cite{chen2013implicit,tadjeran2007second,yang2010numerical})
are oriented to using finite difference approximations for problems
with fractional derivatives in separate directions.
Approximations of fractional derivatives leads to a system of
ordinary differential equations with a filled matrix.
The solution of such problems requires high computational costs
\cite{roop2006computational}. 

To solve problems with fractional powers of elliptic operators, we can apply
finite volume and finite element methods oriented to using arbitrary
domains and irregular computational grids \cite{KnabnerAngermann2003,QuarteroniValli1994}.
The computational realization is associated with the implementation of the matrix function-vector multiplication, i.e., $\varPhi(A) b$. For example, considering the backward Euler scheme,
we have $\varPhi(z) = (1 + \tau z^\alpha)^{-1}$, where $\tau$ is a time step.
To evaluate $\varPhi(A) b$, different approaches \cite{higham2008functions} are available.
Problems of using Krylov subspace methods with the Lanczos approximation
when solving systems of linear equations associated with
the fractional elliptic equations are discussed in \cite{ilic2009numerical}.
A comparative analysis of the contour integral method, the extended Krylov subspace method, and the preassigned poles and interpolation nodes method for solving
space-fractional reaction-diffusion equations is presented in \cite{burrage2012efficient}.
The simplest variant is associated with the explicit construction of the solution using the known
eigenvalues and eigenfunctions of the elliptic operator with diagonalization of the corresponding matrix
\cite{bueno2012fourier,ilic2005numerical,ilic2006numerical}. 
Unfortunately, all these approaches demonstrates too high computational complexity for multidimensional problems.

In the recent paper \cite{vabishchevich2014numerical},
we have proposed a computational algorithm for solving
an equation for fractional powers of elliptic operators on the basis of
a transition to a pseudo-parabolic equation.
For the auxiliary Cauchy problem, the standard two-level schemes are applied.
The computational algorithm is simple for practical use, robust, and applicable to solving
a wide class of problems. A small number of time steps is required to find a solution.
Here this computational algorithm for solving equations with fractional powers of operators
is extended to transient problems. To solve numerically the problem, we construct a special
two-level regularized difference scheme, which is unconditionally stable.

The paper is organized as follows.
The formulation of an unsteady problem containing a fractional power 
of an elliptic operator is given in Section 2. 
Finite element approximation in space is discussed in Section 3. 
In Section 4, we construct a regularized difference scheme 
and investigate its stability. The computational algorithm for solving the equation with a fractional 
power of an operator based on the Cauchy problem for a pseudo-parabolic 
equation is proposed in Section 5.
The results of numerical experiments are described in Section 6.

\section{Problem formulation}\label{sec:2}

\setcounter{section}{2}
\setcounter{equation}{0}\setcounter{theorem}{0}

In a bounded polygonal domain $\Omega \subset R^m$, $m=1,2,3$ with the Lipschitz continuous boundary $\partial\Omega$,
we search the solution for the problem with a fractional power of an elliptic operator.
Define the elliptic operator as
\begin{equation}\label{2.1}
  \mathcal{A}  u = - {\rm div}  k({\bm x}) {\rm grad} \, u + c({\bm x}) u
\end{equation} 
with coefficients $0 < k_1 \leq k({\bm x}) \leq k_2$, $c({\bm x}) \geq 0$.
The operator $\mathcal{A}$ is defined on the set of functions $u({\bm x})$ that satisfy
on the boundary $\partial\Omega$ the following conditions:
\begin{equation}\label{2.2}
  k({\bm x}) \frac{\partial u }{\partial n } + \mu ({\bm x}) u = 0,
  \quad {\bm x} \in \partial \Omega ,
\end{equation} 
where $\mu ({\bm x}) \geq \mu_1 > 0, \  {\bm x} \in \partial \Omega$.

In the Hilbert space $H = L_2(\Omega)$, we define 
the scalar product and norm in the standard way:
\[
  <u,v> = \int_{\Omega} u({\bm x}) v({\bm x}) d{\bm x},
  \quad \|u\| = <u,u>^{1/2} .
\] 
In the spectral problem
\[
 \mathcal{A}  \varphi_k = \lambda_k \varphi_k, 
 \quad \bm x \in \Omega , 
\] 
\[
  k({\bm x}) \frac{\partial  \varphi_k}{\partial n } + \mu ({\bm x}) \varphi_k = 0,
  \quad {\bm x} \in \partial \Omega , 
\] 
we have 
\[
 \lambda_1 \leq \lambda_2 \leq ... ,
\] 
and the eigenfunctions  $ \varphi_k, \ \|\varphi_k\| = 1, \ k = 1,2, ...  $ form a basis in $L_2(\Omega)$. Therefore, 
\[
 u = \sum_{k=1}^{\infty} (u,\varphi_k) \varphi_k .
\] 
Let the operator $\mathcal{A}$ be defined in the following domain:
\[
 D(\mathcal{A} ) = \{ u \ | \ u(x) \in L_2(\Omega), \ \sum_{k=0}^{\infty} | (u,\varphi_k) |^2 \lambda_k < \infty \} .
\] 
Under these conditions  $\mathcal{A} : L_2(\Omega) \rightarrow L_2(\Omega)$ and
the operator $\mathcal{A}$ is self-adjoint and positive defined: 
\begin{equation}\label{2.3}
  \mathcal{A}  = \mathcal{A} ^* \geq \delta I ,
  \quad \delta > 0 ,    
\end{equation} 
where $I$ is the identity operator in $H$.
For $\delta$, we have $\delta = \lambda_1$.
In applications, the value of $\lambda_1$ is unknown (the spectral problem must be solved).
Therefore, we assume that $\delta \leq \lambda_1$ in (\ref{2.3}).
Let us assume for the fractional power of the  operator $A$
\[
 \mathcal{A} ^\alpha u =  \sum_{k=0}^{\infty} (u,\varphi_k) \lambda_k^\alpha  \varphi_k .
\] 
More general and mathematically complete definition of fractional powers of elliptic operators 
is given in \cite{yagi2009abstract}. 

We seek the solution of the Cauchy problem for the evolutionary 
first-order equation with the fractional power of the operator $\mathcal{A}$. 
The solution $u(\bm x,t)$ satisfies the equation
\begin{equation}\label{2.4}
  \frac{d u}{d t} + \mathcal{A}^\alpha u = f(t),
  \quad 0 < t \leq T,  
\end{equation} 
and the initial condition
\begin{equation}\label{2.5}
  u(0) = u_0,
\end{equation} 
under the restriction $0 < \alpha < 1$. 

\section{Discretization in space}\label{sec:3}

\setcounter{section}{3}
\setcounter{equation}{0}\setcounter{theorem}{0}

To solve numerically the problem (\ref{2.4}), (\ref{2.5}), we employ finite element 
approximations in space \cite{brenner2008mathematical,Thomee2006}. 
For (\ref{2.1}) and (\ref{2.2}), we define the bilinear form
\[
 a(u,v) = \int_{\Omega } \left ( k \, {\rm grad} \, u \, {\rm grad} \, v + c \, u v \right )  d {\bm x} +
 \int_{\partial \Omega } \mu \, u v d {\bm x} .
\] 
By (\ref{2.3}), we have
\[
a(u,u) \geq \delta \|u\|^2 .  
\]
Define a subspace of finite elements $V^h \subset H^1(\Omega)$.
Let $\bm x_i, \ i = 1,2, ..., M_h$ be triangulation points for the domain $\Omega$.
Define pyramid function $\chi_i(\bm x) \subset V^h, \ i = 1,2, ..., M_h$, where
\[
 \chi_i(\bm x_j) = \left \{
 \begin{array}{ll} 
 1, & \mathrm{if~}  i = j, \\
 0, & \mathrm{if~}  i \neq  j .
 \end{array}
 \right . 
\] 
For $v \in V_h$, we have
\[
 v(\bm x) = \sum_{i=i}^{M_h} v_i \chi_i(\bm x),
\] 
where $v_i = v(\bm x_i), \ i = 1,2, ..., M_h$.

We define the discrete elliptic operator $A$ as
\[
 a(u,v) = \ <Au, v>,
 \quad \forall \ u,v \subset V^h , 
\] 
where, similarly to (\ref{2.3}),  
\begin{equation}\label{3.1}
 A = A^* \geq \delta I ,
 \quad \delta > 0 .  
\end{equation} 
For the problem (\ref{2.4}), (\ref{2.5}), we put into the correspondence the operator equation 
for $w(t) \subset V^h$:
\begin{equation}\label{3.2}
 \frac{d w}{d t} + A^\alpha w = \psi(t), 
 \quad 0 < t \leq T, 
\end{equation} 
\begin{equation}\label{3.3}
 w(0) = w_0, 
\end{equation} 
where $\psi(t) = P f(t)$, $w_0 = P u_0$ 
with $P$ denoting $L_2$-projection onto $V^h$.

Now we will obtain an elementary a priori estimate for the solution of (\ref{3.2}), (\ref{3.3})
assuming that the solution of the problem, coefficients of the elliptic operator, the right-hand side and initial conditions are sufficiently smooth.

Let us multiply equation (\ref{3.2}) by $w$ and integrate it over the domain $\Omega$:
\[
 \left <\frac{d w}{d t}, w \right > + < A^\alpha w, w > \, = \,
 < \psi, w > .
\]
In view of the self-adjointness and positive definiteness of the operator $A^\alpha$,  
the right-hand side can be evaluated by the inequality
\[
 < \psi, w > \, \leq \, <A^\alpha w, w > + \frac{1}{4}  <A^{-\alpha} \psi, \psi  > .
\] 
By virtue of this, we have
\[
 \frac{d}{d t} \|w\|^2 \leq \frac{1}{2} \|\psi \|^2_{A^{-\alpha}} ,
\] 
where $\|\psi \|_{A^{-\alpha}} = <A^{-\alpha} \psi, \psi >^{1/2}$. 
The latter inequality leads us to the desired a priori estimate:
\begin{equation}\label{3.4}
 \|w(t)\|^2 \leq \|w_0\|^2 + \frac{1}{2} \int_{0}^{t}\|\psi(\theta) \|^2_{A^{-\alpha}} d \theta .
\end{equation} 
Taking into account (\ref{3.1}),  the estimate (\ref{3.4}) can be simplified:
\begin{equation}\label{3.5}
 \|w(t)\|^2 \leq \|w_0\|^2 + \frac{1}{2} \delta^{-\alpha} \int_{0}^{t}\|\psi(\theta) \|^2 d \theta .
\end{equation}
We will focus on the estimates (\ref{3.4}), (\ref{3.5}) 
for the stability of the solution with respect to the initial data and the right-hand side 
in constructing discrete analogs of the problem (\ref{3.2}), (\ref{3.3}).

 \section{Regularized scheme}\label{sec:4}

\setcounter{section}{4}
\setcounter{equation}{0}\setcounter{theorem}{0}

To solve numerically the problem (\ref{3.2}), (\ref{3.3}), we use the simplest implicit two-level scheme.
Let $\tau$ be a step of a uniform grid in time such that $w^n = w(t^n), \ t^n = n \tau$,
$n = 0,1, ..., N, \ N\tau = T$.
It seems reasonable to begin with the simplest explicit scheme
\begin{equation}\label{4.1}
 \frac{w^{n+1} - w^{n}}{\tau } + A^\alpha w^{n} = \psi^{n},
 \quad n = 0,1, ..., N-1,
\end{equation} 
\begin{equation}\label{4.2}
 w^0 = w_0 .
\end{equation}
Advantages and disadvantages of explicit schemes for the standard parabolic problem ($\alpha = 1$)
are well-known, i.e., these are a simple computational implementation and a time step restriction
(see, e.g., \cite{Samarskii1989,SamarskiiMatusVabischevich2002}).
In our case ($\alpha \neq 1$), the main drawback (conditional stability) 
remains, whereas the advantage in terms of implementation simplicity does not exist. 
The approximate solution at a new time level is determined via (\ref{4.1}) as
\begin{equation}\label{4.3}
 w^{n+1} = w^{n} - \tau A^\alpha w^{n} + \tau \psi^{n} . 
\end{equation} 
Thus, we must calculate $A^\alpha w^{n}$.
In view of these problems, considering the scheme (\ref{4.1}), it is more correct to speak of 
the scheme with the explicit approximations in time in contrast to the standard fully explicit scheme.

Let us approximate equation (\ref{3.2}) by the backward Euler scheme:
\begin{equation}\label{4.4}
 \frac{w^{n+1} - w^{n}}{\tau } + A^\alpha w^{n+1} = \psi^{n+1},
 \quad n = 0,1, ..., N-1.
\end{equation}
The main advantage of the implicit scheme (\ref{4.4}) in comparison with (\ref{4.1})
is its absolute stability. Let us derive for this scheme  the corresponding estimate for stability.

Multiplying equation (\ref{4.4}) scalarly by $\tau w^{n+1}$, we obtain
\begin{equation}\label{4.5}
\begin{split}
 < w^{n+1}, w^{n+1}> \ & + \ \tau < A^\alpha w^{n+1}, w^{n+1}> \  \\
 & = \ < w^{n}, w^{n+1}> + \tau < \psi^{n+1}, w^{n+1}> .
\end{split} 
\end{equation}
The terms on the right side of (\ref{4.5}) are estimated using the inequalities:
\[
 < w^{n}, w^{n+1}> \ \leq \frac{1}{2} < w^{n+1}, w^{n+1}> + \frac{1}{2} < w^{n}, w^{n}> ,
\] 
\[
 < \psi^{n+1}, w^{n+1}> \ \leq \ < A^\alpha w^{n+1}, w^{n+1}> + \frac{1}{4}  < A^{-\alpha} \psi^{n+1}, \psi^{n+1}> .
\]
The substitution into (\ref{4.5}) leads to the following level-wise estimate:
\[
 \|w^{n+1}\|^2 \leq \|w^{n}\|^2 + \frac{1}{2} \tau \|\psi^{n+1}\|^2_{A^{-\alpha}} .
\]
This implies the desired estimate for stability:
\begin{equation}\label{4.6}
 \|w^{n+1}\|^2 \leq \|w_0\|^2 + \frac{1}{2} \sum_{k=0}^{n}\tau \|\psi^{k+1}\|^2_{A^{-\alpha}} ,
\end{equation} 
which is a discrete analog of the estimate (\ref{3.4}). 
Similarly to (\ref{3.5}), in view of  (\ref{3.1}), from (\ref{4.6}), we get
\begin{equation}\label{4.7}
 \|w^{n+1}\|^2 \leq \|w_0\|^2 + \frac{1}{2} \delta^{-\alpha} \sum_{k=0}^{n}\tau \|\psi^{k+1}\|^2 .
\end{equation}
To obtain the solution at the new time level, it is necessary to solve the problem
\[
 (I +  \tau A^\alpha) w^{n+1} = w^{n} + \tau \psi^{n} .
\]
In our case, we must calculate the values of $\varPhi(A) b$ for $\varPhi(z) = (1+ \tau z^\alpha)^{-1}$. 

A more complicated situation arises in the implementation of the Crank-Nicolson scheme:
\[
 \frac{w^{n+1} - w^{n}}{\tau } + A^\alpha  \frac{w^{n+1} + w^{n}}{2} = \psi^{n+1/2},
 \quad n = 0,1, ..., N-1.
\]
In this case, we have 
\[
 \left (I +  \frac{\tau}{2}  A^\alpha \right ) w^{n+1} = w^{n} - \frac{\tau}{2}  A^\alpha   w^{n}
 + \tau \psi^{n+1/2} ,
\]
i.e., we need to evaluate both $\varPhi(z) = (1+ 0.5 \tau z^\alpha)^{-1}$ and $\varPhi(z) = z^\alpha$.

The numerical implementation of the above-mentioned approximations in time for the standard parabolic problems ($\alpha = 1$ in (\ref{3.2})) is based on calculating the values of $\varPhi(A) b$ 
for $\varPhi(z) = (1+ \sigma \tau z)^{-1}, \ \sigma =0.5, 1$ and $\varPhi(z) = z$. 
For problems with fractional powers of elliptic operators, 
we apply the approach proposed early in our paper \cite{vabishchevich2014numerical}.
It is based on the computation of $\varPhi(A) b$ for $\varPhi(z) = z^{-\beta}, \ 0 < \beta < 1$.

For the explicit approximation in time, we rewrite (\ref{4.3}) in the form
\[
 w^{n+1} = w^{n} - \tau A A^{-\beta} w^{n} + \tau \psi^{n} ,
 \quad \beta = 1 - \alpha .   
\]
Therefore, the computational implementation is based on the evaluation of 
$\varPhi(A) b$ for $\varPhi(z) = z^{-\beta}$ and $\varPhi(z) = z$. 
A similar approach is not valid for the  backward Euler scheme (\ref{4.2}), (\ref{4.4}) 
and moreover for the Crank-Nicolson scheme. To construct
a more appropriate from a computational point of view 
approximations in time for the Cauchy problem (\ref{3.2}), (\ref{3.3}), we apply 
the principle of regularization for operator-difference schemes proposed by A.A. Samarskii \cite{Samarskii1989}.

For a regularizing operator $R = R^* > 0$, the simplest regularized scheme 
for solving (\ref{3.2}), (\ref{3.3}) has the form (see, e.g., \cite{Vabishchevich2014}):
\begin{equation}\label{4.8}
 (I + \tau R) \frac{w^{n+1} - w^{n}}{\tau } + A^\alpha w^{n} = \psi^{n+1},
 \quad n = 0,1, ..., N-1.
\end{equation}
Now we will derive the stability conditions for the regularized scheme (\ref{4.2}), (\ref{4.8})  
and after that we will select the appropriate regularizing operator $R$ itself.

Rewrite equation (\ref{4.8}) in the form
\[
 \left (I + \tau \left ( R - \frac{1}{2}A^\alpha \right ) \right ) \frac{w^{n+1} - w^{n}}{\tau } 
 + A^\alpha \frac{w^{n+1} + w^{n}}{2} = \psi^{n+1} .
\] 
Multiplying it scalarly by $\tau (w^{n+1} + w^{n})$, we get
\[
\begin{split}
 < D w^{n+1}, w^{n+1}> - < D w^{n}, w^{n}> & +  \frac{\tau }{2} <A^\alpha (w^{n+1} + w^{n}), w^{n+1} + w^{n} > \\
 & = \tau <\psi^{n+1}, w^{n+1} + w^{n} > ,
\end{split}
\] 
where
\begin{equation}\label{4.9}
 D = I + \tau \left ( R - \frac{1}{2}A^\alpha \right ) .
\end{equation} 
For
\begin{equation}\label{4.10}
 R \geq \frac{1}{2}A^\alpha 
\end{equation} 
we have $D = D^* \geq I$, and $D = I + \mathit{O}(\tau)$.
Under these conditions, we obtain the inequality
\[
 \|w^{n+1}\|_D^2 \leq \|w^{n}\|_D^2 + \frac{1}{2} \tau \|\psi^{n+1}\|^2_{A^{-\alpha}} .
\]
Thus, for the regularized difference scheme (\ref{4.2}), (\ref{4.8}), 
under the condition (\ref{4.10}), the following estimate for stability with respect to the initial data 
and the right-hand side holds:
\begin{equation}\label{4.11}
 \|w^{n+1}\|_D^2 \leq \|w_0\|_D^2 + \frac{1}{2} \sum_{k=0}^{n}\tau \|\psi^{k+1}\|^2_{A^{-\alpha}} .
\end{equation} 

To select an appropriate regularizing operator $R$, we should take into account two conditions, i.e.,
first, to satisfy the inequality (\ref{4.10}), and secondly, to simplify calculations. 
Our choice is based on the inequality
\begin{equation}\label{4.12}
 A^\alpha \leq \alpha A + (1-\alpha) I ,
\end{equation}
which is the simplest version of Young's inequality for positive operators 
(see, e.g., \cite{FangDu}). In the scheme (\ref{4.8}), we put
\begin{equation}\label{4.13}
 R = \sigma (\alpha A + (1-\alpha) I) .
\end{equation}
For $\sigma \geq 0.5$, in view of (\ref{4.12}), the inequality (\ref{4.10}) holds.

The result of our analysis is the following statement.

\begin{theorem}\label{Th1}
The regularized scheme (\ref{4.2}), (\ref{4.8})  with the regularizer $R$
selected according to (\ref{4.13}) is unconditionally stable for $\sigma \geq 0.5$.
The approximate solution satisfies the a priori estimate (\ref{4.9}), (\ref{4.11}). 
\end{theorem}

The transition to a new time level is performed via the formula
\[
\begin{split}
 (I +  & \sigma \tau (\alpha A + (1-\alpha) I ) w^{n+1} =  
 (I +  \sigma \tau (\alpha A + (1-\alpha) I ) w^{n}  \\
 & - \tau A A^{-\beta} w^{n} + \tau \psi^{n+1} ,
 \quad \beta = 1 - \alpha . 
\end{split}
\]
Therefore,  it is necessary to calculate the values 
$\varPhi(A) b$ for $\varPhi(z) = (1+\tau\widetilde{\sigma} z)^{-1}$ 
and $\varPhi(z) = z^{-\beta}, \ 0 < \beta < 1$.

\section{Calculation of the operator with the fractional power}\label{sec:5}

\setcounter{section}{5}
\setcounter{equation}{0}\setcounter{theorem}{0}

The main peculiarity of solving the  Cauchy problem  (\ref{3.2}), (\ref{3.3}) 
is the necessity to evaluate values
\[
 g^n = A^{-\beta} w^{n}, 
 \quad n = 0,1, ..., N-1,
 \quad  0 < \beta = 1 - \alpha  < 1 .
\]
The computational algorithm is based on the consideration of the auxiliary 
Cauchy problem \cite{vabishchevich2014numerical}.

Assume that
\[
 y(s) = \delta^{\beta} (s (A - \delta I) + \delta I)^{-\beta} y(0) ,
\]
then for the determination of $g^n$, we can put
\begin{equation}\label{5.1}
 g^n = y(1),
 \quad y(0) = \delta^{-\beta}  w^{n}.  
\end{equation} 
The function $y(s)$ satisfies the evolutionary equation
\begin{equation}\label{5.2}
  (s G + \delta I) \frac{d y}{d s} + \beta G y = 0 ,
\end{equation} 
where $G = A - \delta I \geq  0$.
Thus, the calculation of $A^{-\beta} w^{n}$  is based on the solution of the Cauchy problem (\ref{5.1}), (\ref{5.2}) within  the unit interval for the pseudo-parabolic equation.

To solve numerically the problem (\ref{5.1}), (\ref{5.2}),
we use the simplest implicit two-level scheme.
Let $\eta$ be a step of a uniform grid in time such that $y_k = y(s_k), \ s_k = k \eta$,
$k = 0,1, ..., K, \ K\eta = 1$.
Let us approximate equation (\ref{5.2})  by the backward Euler scheme
\begin{equation}\label{5.3}
  (s_{k+1} G + \delta I) \frac{ y_{k+1} - y_{k}}{\eta } + 
  \beta  G  y_{k+1}  = 0,  
  \quad k = 0,1, ..., K-1,
\end{equation}
\begin{equation}\label{5.4}
  y_0 = \delta^{-\beta } w^{n} . 
\end{equation}
For the Crank-Nicolson scheme, we have
\begin{equation}\label{5.5}
  (s_{k+1/2} G + \delta I) \frac{ y_{k+1} - y_{k}}{\eta } + 
  \beta  G  \frac{ y_{k+1} + y_{k}}{2}  = 0,  
  \quad k = 0,1, ..., K-1 .
\end{equation}
The difference scheme (\ref{5.4}), (\ref{5.5}) approximates the problem 
(\ref{5.1}), (\ref{5.2})
with the second order by $\eta $, whereas for scheme (\ref{5.3}), (\ref{5.4})  we have only the first order.

The above two-level schemes are unconditionally stable. 
The corresponding level-wise estimate has the form
\begin{equation}\label{5.6}
 \|y_{k+1}\| \leq \|y_{k}\| ,
 \quad k = 0,1, ..., K-1 . 
\end{equation}
To prove (\ref{5.6}) (see \cite{vabishchevich2014numerical}), it is sufficient to multiply scalarly
equation (\ref{5.3}) by $y_{k+1}$ and equation (\ref{5.5}) by $y_{k+1}+y_{k}$. 
Taking into account (\ref{5.4}), from (\ref{5.6}), we obtain
\begin{equation}\label{5.7}
 \|y_K\| \leq \delta^{-\beta } \|w^{n}\| . 
\end{equation} 

The solution of the Cauchy problem (\ref{5.3}), (\ref{5.4}) (or (\ref{5.4}), (\ref{5.5})) 
may be written in the form
\begin{equation}\label{5.8}
 y_K = Q_K(A) w^{n} .
\end{equation}
For instance, for the backward Euler scheme  (\ref{5.3}), (\ref{5.4}), we have
\begin{equation}\label{5.9}
 Q_K(A) = \prod_{k = 1}^{K} (k \eta G + \delta I +  \eta  \beta  G)^{-1} (k \eta G + \delta I) .
\end{equation} 
As for the Crank-Nicolson scheme (\ref{5.3}), (\ref{5.4}), we obtain the representation 
\begin{equation}\label{5.10}
 Q_K(A) = \prod_{k = 1}^{K} (k \eta G + \delta I +  0.5 \eta   \beta  G)^{-1} (k \eta G + \delta I - 0.5 \eta   \beta  G) .
\end{equation} 

The computational implementation of the regularized scheme (\ref{4.2}), (\ref{4.8}), (\ref{4.13}) 
based on solving the auxillary evolutionary problems (\ref{5.1}), (\ref{5.2}) 
corresponds to the new scheme
\begin{equation}\label{5.11}
 (I + \tau \sigma (\alpha A + (1-\alpha) I)) \frac{w^{n+1} - w^{n}}{\tau } + A Q_K(A) w^{n} = \psi^{n+1} .
\end{equation} 
 
 \begin{theorem}\label{Th2}
The scheme (\ref{4.2}), (\ref{5.9}), (\ref{5.11}) is unconditionally stable for
\begin{equation}\label{5.12}
 \sigma \geq \frac{1}{2} + \frac{1-\alpha }{2\alpha \delta } .
\end{equation}
Moreover, the solution satisfies the stability estimate  
\begin{equation}\label{5.13}
 \|w^{n+1}\|_D \leq \|w_0\|_D + \sum_{k=0}^{n}\tau \|\psi^{k+1}\|_{D^{-1}} ,
\end{equation} 
where $D \geq I$ and
\begin{equation}\label{5.14}
 D = I + \tau \left (\sigma (\alpha A + (1-\alpha) I) - \frac{1}{2}A Q_K(A) \right ) .
\end{equation} 
 \end{theorem}

 \proof 
First, we will show that under the above restrictions on $\sigma$ (\ref{5.12}), 
we have that the operator $D \geq I$.  According to (\ref{5.7}), for  (\ref{5.9}), we have
\begin{equation}\label{5.15}
 0 < Q_K(A) \leq  \delta^{-\beta } I . 
\end{equation} 
By (\ref{3.1}), (\ref{5.15}) and Youngs inequality, we obtain
\[
\begin{split}
 \sigma (\alpha A + (1-\alpha) I) & - \frac{1}{2}A Q_K(A) \geq   \sigma (\alpha A + (1-\alpha) I) - \frac{1}{2} \delta^{-\beta } A \\
 & \geq (2 \delta)^{-1} (2 \sigma \alpha \delta - \delta^{\alpha}) A  \\
 & \geq (2 \delta)^{-1} ((2 \sigma - 1) \alpha \delta - (1-\alpha))  A \geq 0 .
\end{split}
\] 
Next, rewrite the scheme (\ref{5.11}) in the form
\[
 D \frac{w^{n+1} - w^{n}}{\tau } +  A Q_K(A) \frac{w^{n+1} + w^{n}}{2}  = \psi^{n+1} .
\]
Multiplying this equation scalarly by $\tau (w^{n+1} + w^{n})$, in view of the left inequality (\ref{5.15}), we arrive at
\[
 \|w^{n+1}\|_D^2 - \|w^{n}\|_D^2 \leq \tau <\psi^{n+1}, w^{n+1} + w^{n} > .
\] 
Taking into account
\[
 \|w^{n+1}\|_D^2 - \|w^{n}\|_D^2 = (\|w^{n+1}\|_D - \|w^{n}\|_D) (\|w^{n+1}\|_D + \|w^{n}\|_D),
\] 
\[
 <\psi^{n+1}, w^{n+1} + w^{n} > \ \leq \|\psi^{n+1}\|_{D^{-1}} (\|w^{n+1}\|_D + \|w^{n}\|_D),
\] 
we get the estimate
\[
 \|w^{n+1}\|_D \leq \|w^{n}\|_D + \tau  \|\psi^{n+1}\|_{D^{-1}} .
\]
These inequalities prove the estimate (\ref{5.13}), (\ref{5.14}) for stability 
of the difference scheme with respect on the initial data and the right-hand side.
 \proofend 

Thus, the stability of the difference scheme with an approximate calculation of  $A^{-\beta} w^{n}$
via the backward Euler scheme (\ref{5.3}), (\ref{5.4}) is proved under  more
strong restrictions on the weight parameter $\sigma$. 
In the original regularized scheme (see Theorem~\ref{Th1}), 
it was enough to take $\sigma \geq 0.5$, whereas here we have (\ref{5.12}).

As for the Crank-Nicolson scheme (\ref{5.4}), (\ref{5.5}), for the approximate evaluation 
of $A^{-\beta} w^{n}$ in the scheme (\ref{5.11}),  the operator $Q_K(A)$  is determined according
to (\ref{5.10}).
This operator is no longer positive, i.e., instead of (\ref{5.15}), we have the bilateral inequality
\[
 - \delta^{-\beta } I  \leq  Q_K(A) \leq  \delta^{-\beta } I . 
\]
In this case, it is no possible to establish the unconditional stability of the scheme
(\ref{4.2}), (\ref{5.9}).

\section{Numerical experiments}\label{sec:6}

\setcounter{section}{6}
\setcounter{equation}{0}\setcounter{theorem}{0} 

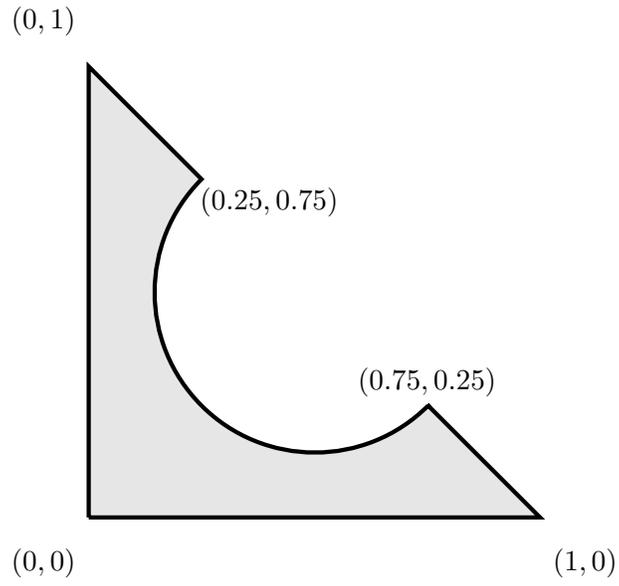
\begin{figure}[htp]
  \begin{center}
  \begin{tikzpicture}[scale = 0.6]
     \draw [ultra thick, fill=gray!20] (0,0) -- (0,10) -- (2.5,7.5)  arc [radius=3.55, start angle=135, end angle= 315] -- (7.5,2.5) -- (10,0) --  (0,0);
     \draw(-1,-1) node {$(0,0)$};
     \draw(11,-1) node {$(1,0)$};
     \draw(4,7) node {$(0.25,0.75)$};
     \draw(-1,11) node {$(0,1)$};
     \draw(7.5,3) node {$(0.75,0.25)$};
  \end{tikzpicture}
  \caption{Computational domain $\Omega$}
  \label{f-1}
  \end{center}
\end{figure}

\begin{figure}[!h]
  \begin{center}
    \includegraphics[scale=0.2] {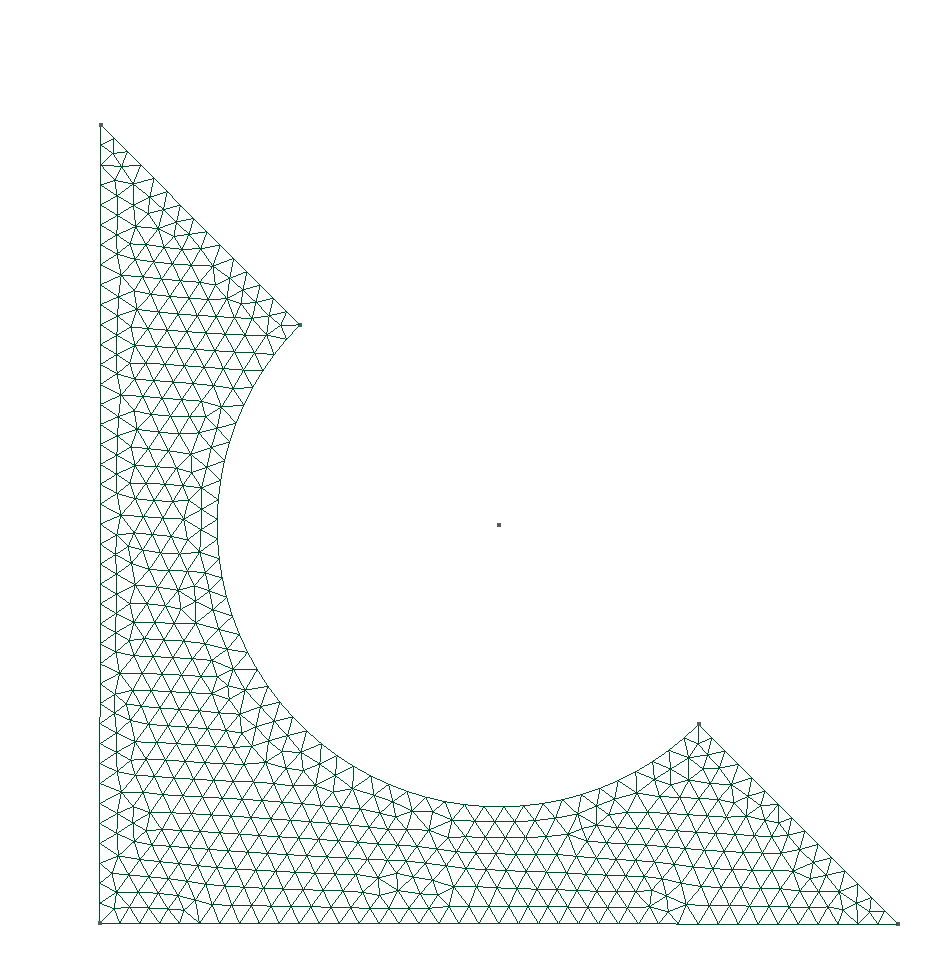}
	\caption{Medium grid 2: 679 nodes, 1201 triangles}
	\label{f-2}
  \end{center}
\end{figure}

Capabilities of the proposed method are illustrated by solving a model two-dimensional  
problem. The computational domain is shown in Fig.~\ref{f-1}. 
Triangulation is performed to discretize this domain. 
Calculations are performed using the  coarse (grid 1: 198 nodes, 315 triangles), 
medium (see Fig.~\ref{f-2}) and fine (grid 3: 2470 nodes, 4631 triangles) grids.

The unsteady problem (\ref{2.4}), (\ref{2.5}) is considered for for the elliptic operator (\ref{2.1}), (\ref{2.2}) with constant coefficients:
\[
 k({\bm x}) = 1, 
 \quad c({\bm x}) = 0,
 \quad \mu ({\bm x}) = \mu. 
\]
The right-hand side and the initial condition are given as
\begin{equation}\label{6.1}
 f(\bm x, t) = \frac{2}{1 + \exp(\gamma(x_1-x_2))} ,
 \quad u_0(\bm x) = 0 . 
\end{equation} 
For $\gamma \rightarrow \infty$, the right-hand side becomes discontinuous.
Finite element approximations lead to the Cauchy problem (\ref{3.2}), (\ref{3.3}). 

\clearpage

To estimate the constant $\delta$ in (\ref{3.1}), we solve the spectral problem 
\begin{equation}\label{6.2}
 a(v,v) = \lambda <v,v>,
 \quad v \in V^h,  
\end{equation} 
where $\delta = \lambda_{1}$.
When choosing  piecewise linear finite elements 
($V^h \subset H^1(\Omega)$), the corresponding values of the constant $\delta$ 
for the above-mentioned computational grids and
$\mu = 1, 10, 100$ are presented in Table~\ref{tbl-1}.
The results show that for the evaluation of $\delta$, we  can use
the solution of the incomplete eigenvalue problem obtained on the coarse grid. 
If we use standard algorithms of inverse iteration \cite{bjorck2015numerical,golub2012matrix}, 
then computational costs are not significant.

\begin{table}[!h]
\caption{Constant $\lambda_{1}$}
\begin{center}
\begin{tabular}{|c|c|c|c|} \hline
Grid    & \multicolumn{3}{c|}{ $\mu$ } \\ \cline{2-4}
    & 1 & 10 & 100  \\ \hline
 1   & 11.413872661 & 73.277955924 &  167.20923186 \\
 2   & 11.422023432 & 72.928682512 &  164.36008412  \\
 3   & 11.424101908 & 72.839189190 &  163.65578746 \\ \hline
\end{tabular}
\end{center}
\label{tbl-1}
\end{table}

Here we present some numerical results for the stationary problem
\[
  A^\beta  y = P f. 
\]
Just this problem is solved at each time level in 
the regularization scheme (\ref{4.2}), (\ref{4.8}).  First, we consider the problem with 
the constant righ-hand side ($\gamma = 0$ in (\ref{6.1})) and $\beta = 0.5$. 
Calculations are performed on grid 2 with $\mu =10$.

\begin{figure}[!h]
  \begin{center}
    \includegraphics[width=0.7\linewidth] {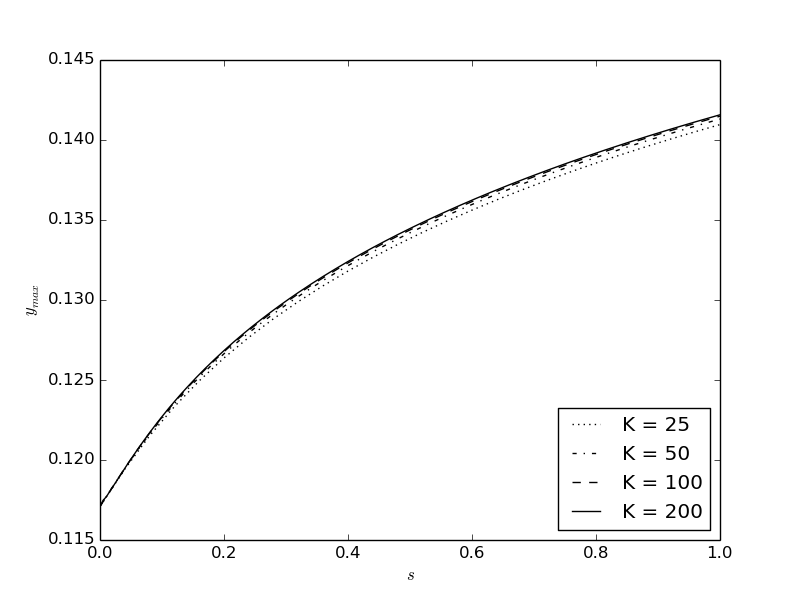}
	\caption{Dynamics of the solution for the Euler scheme ($\delta = \lambda_1$)}
	\label{f-3}
  \end{center}
\end{figure}

\begin{figure}[!h]
  \begin{center}
    \includegraphics[width=0.7\linewidth] {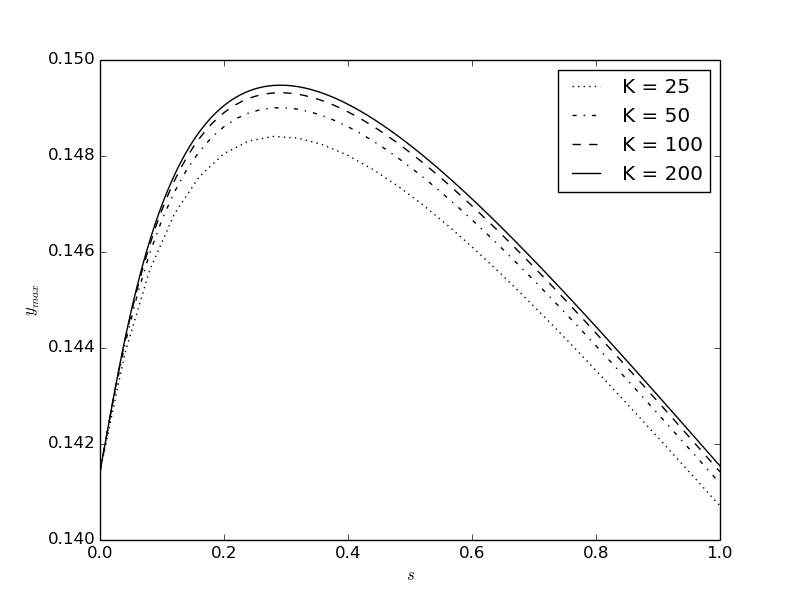}
	\caption{Dynamics of the solution for the Euler scheme ($\delta = 50$)}
	\label{f-4}
  \end{center}
\end{figure}

The most interesting fact for this problem is the dependence on the time step.
Figure~\ref{f-3} presents the dependence of the maximum (over the entire computational domain) value
of the approximate solution on the time step for the backward Euler scheme (\ref{5.3}).
In this case, we used $\delta = \lambda_1 \approx 72.928682512$. 
Figure~\ref{f-4} shows that the parameter $\delta$ demonstrates practically no influence  on the solution.
This calculation was performed with $\delta = 50$.

\begin{figure}[!h]
  \begin{center}
    \includegraphics[width=0.7\linewidth] {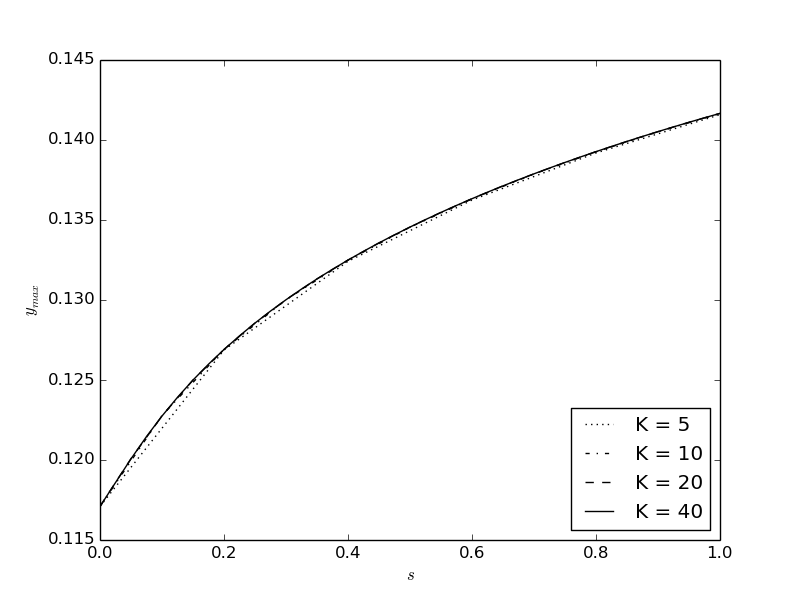}
	\caption{Dynamics of the solution for the Crank-Nicolson scheme ($\delta = \lambda_1$)}
	\label{f-5}
  \end{center}
\end{figure}

\begin{figure}[!h]
  \begin{center}
    \includegraphics[width=0.7\linewidth] {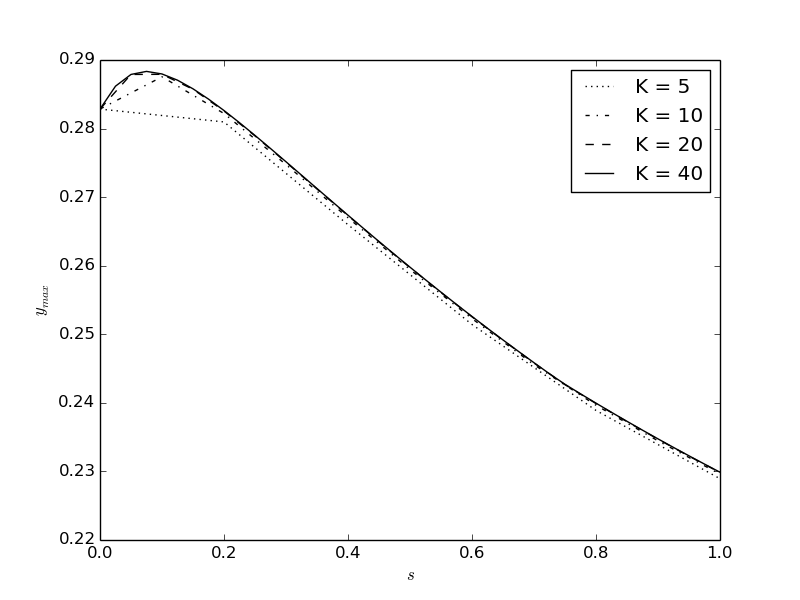}
	\caption{Dynamics of the solution for the Crank-Nicolson scheme ($\delta = 50$)}
	\label{f-6}
  \end{center}
\end{figure}

\begin{figure}[!h]
  \begin{center}
    \includegraphics[width=0.6\linewidth] {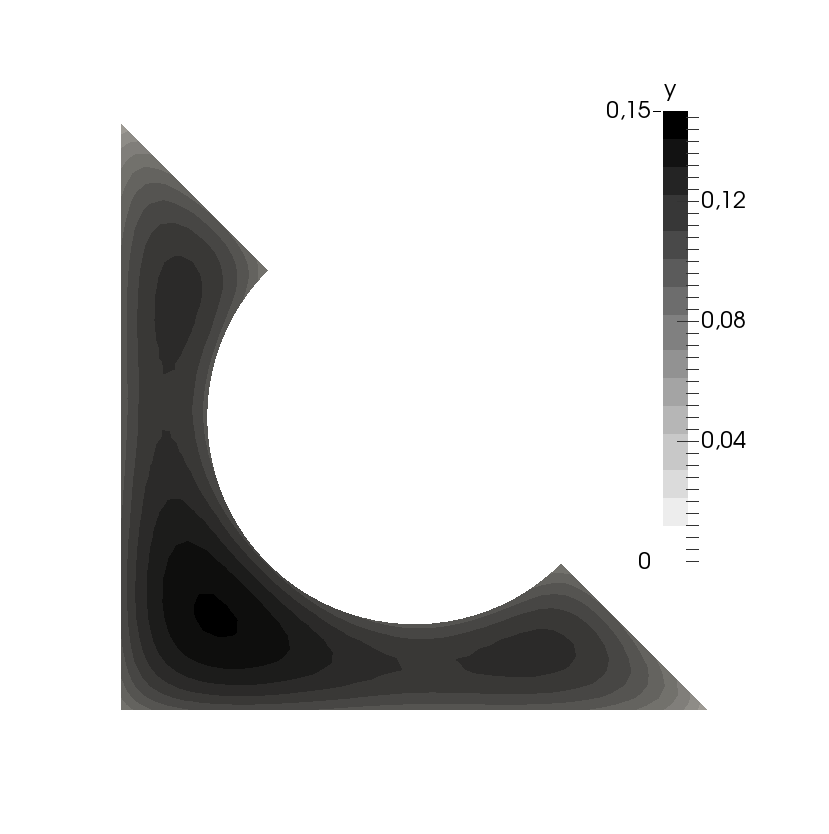}
	\caption{Solution on grid 2 ($y_{max} = 0.141663$)}
	\label{f-7}
  \end{center}
\end{figure}

The convergence of the approximate solution with the first order in time is observed in these figures.
Similar data are depicted in Fig.~\ref{f-5},~\ref{f-6}  for the symmetric scheme, i.e., the Crank-Nicolson scheme (\ref{5.5}). Here we see much more rapid convergence and so
we can obtain acceptable in accuracy results using fairly coarse meshes.
The approximate solution itself is given in Fig.~\ref{f-7}.

\begin{figure}[!h]
  \begin{center}
    \includegraphics[width=0.6\linewidth] {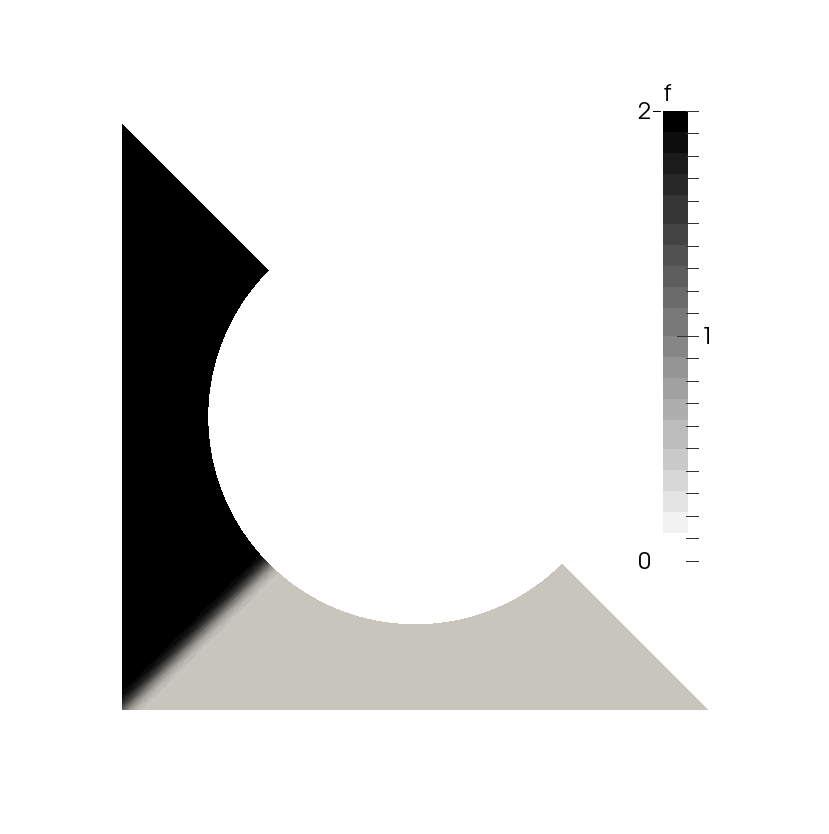}
	\caption{The right-hand side for $\gamma = 100$}
	\label{f-8}
  \end{center}
\end{figure}

\begin{figure}[!h]
  \begin{center}
    \includegraphics[width=0.6\linewidth] {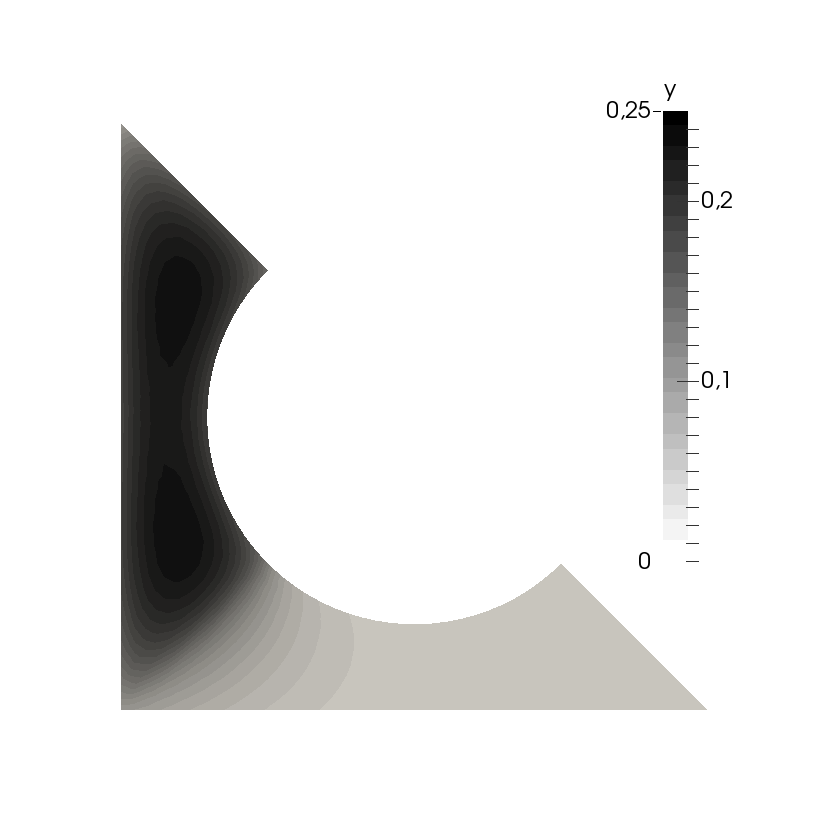}
	\caption{Solution on grid 2 for $\gamma = 100$}
	\label{f-9}
  \end{center}
\end{figure}

\begin{figure}[!h]
  \begin{center}
    \includegraphics[width=0.7\linewidth] {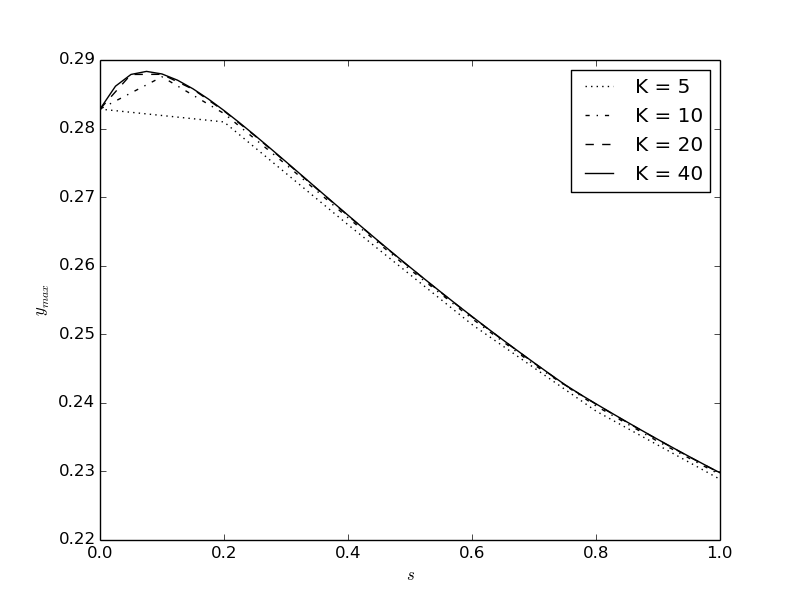}
	\caption{The Crank-Nicolson scheme ($\delta = 50$, $\gamma = 100$)}
	\label{f-10}
  \end{center}
\end{figure}

The effect of the right-hand side is illustrated by the calculations with various values of $\gamma$, which are depicted  in (\ref{6.1}).
For $\gamma = 100$, the right-hand side has the form shown in Fig.~\ref{f-8}.
Figure~\ref{f-9} shows the approximate solution.
The dynamics of the maximum value of the solution is presented in Fig.~\ref{f-10}.
Thus, the calculations demonstrate the high accuracy of the computational algorithm for solving the equation with fractional powers of elliptic operators via the Crank-Nicolson scheme. Moreover, they show
a weak dependence of the accuracy of the approximate solution on the parameter $\delta$ from (\ref{3.1})
as well as on the smoothness of the right-hand side.

\clearpage


\begin{figure}[!h]
  \begin{center}
    \includegraphics[width=0.7\linewidth] {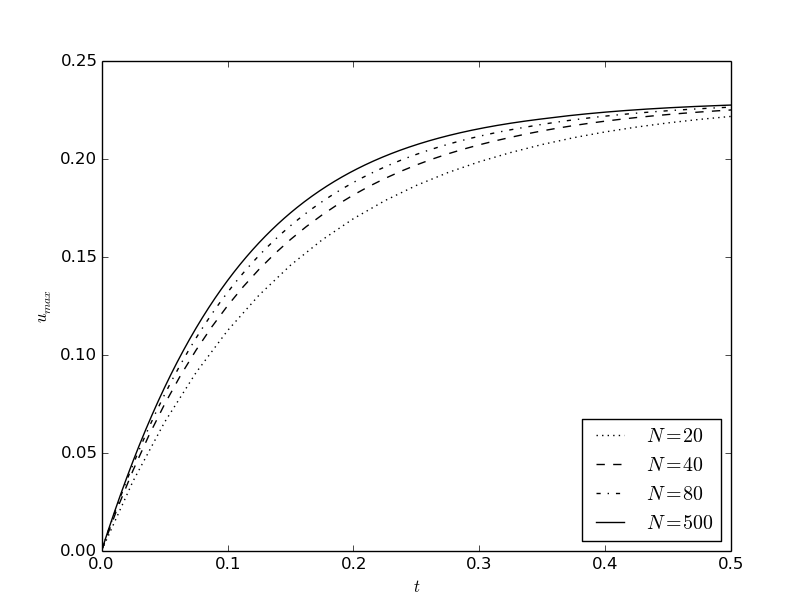}
	\caption{Unsteady problem solution}
	\label{f-11}
  \end{center}
\end{figure}

Now we discuss the numerical results for unsteady problem with $\gamma = 100$  
in (\ref{6.1}) and $\sigma = 0.5$. The achievement of the steady-state solution when $\alpha = 0.5$ is shown in Fig.~\ref{f-11}. We observe the convergence of the approximate solution with the first order by $\tau$.
The calculations were performed using the Crank-Nicolson scheme (\ref{5.3}), (\ref{5.4})   
with $K=10$.

\begin{figure}[!h]
  \begin{center}
    \includegraphics[width=0.7\linewidth] {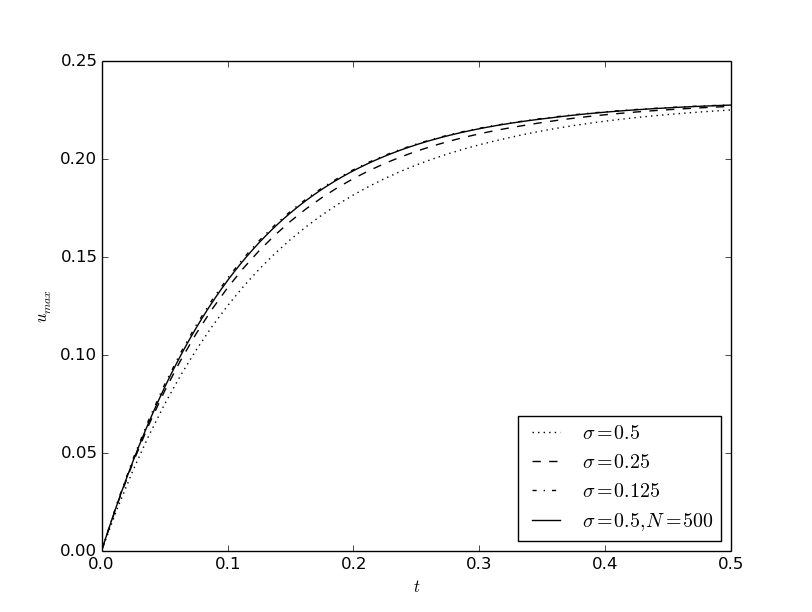}
	\caption{Predictions with various $\sigma$}
	\label{f-12}
  \end{center}
\end{figure}

The condition $\sigma \geq 0.5$ is sufficient for the unconditional stability 
of the  regularized scheme (\ref{4.2}), (\ref{4.8}), (\ref{4.13}).
Trancation error increases when the value of the parameter $\sigma$ becomes higher. 
Data depicted in Fig.~\ref{f-12} demonstrate the effect of $\sigma$.
The problem is solved on the time grid with  $N = 40$, and for a comparison, data 
predicted on the fine grid $N = 500$ with $\sigma = 0.5$ is shown, too. 
Note that in this example, the instability takes place only if 
$\sigma \geq \sigma^* \approx 0.01$. Thus, the sufficient condition for stability 
seems to be essentially exaggerated.

\begin{figure}[!h]
  \begin{center}
    \includegraphics[width=0.7\linewidth] {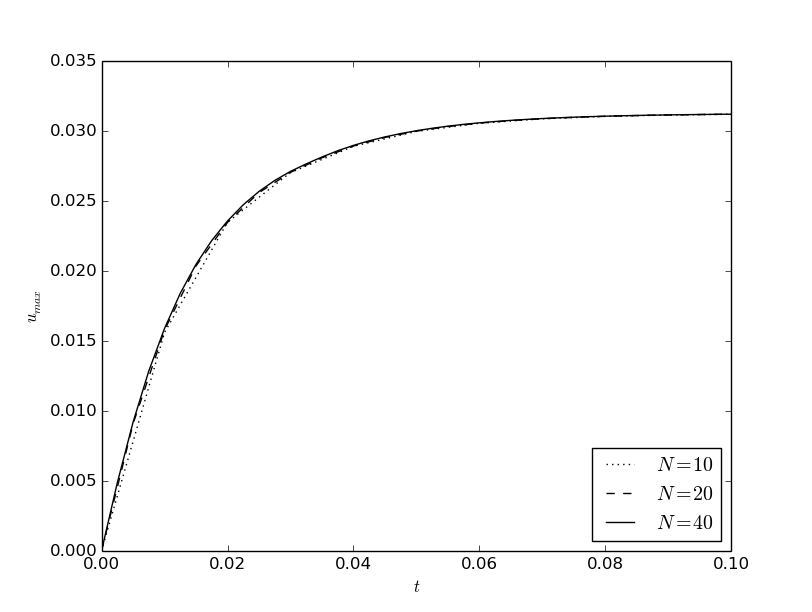}
	\caption{Solution of the unsteady problem for $\alpha = 0.95$}
	\label{f-13}
  \end{center}
\end{figure}

In problems, which are closer to the standard problems of unsteady diffusion 
($\alpha \rightarrow 1$), restrictions on  $\sigma$ seems to be close to optimal ones. 
For example, Figure~\ref{f-13} presents the calculations with $\alpha = 0.95$
and t$\sigma = 0.5$.
In this case, accuracy is much higher, and $\sigma^* \approx 0.35$ for the grid of $N = 40, T = 0.1$. 

\section*{Acknowledgements}

This work was supported by the Russian Foundation for Basic Research  (project 14-01-00785).

\clearpage



 \bigskip \smallskip

 \it

 \noindent
$^1$ Nuclear Safety Institute \\
Russian Academy of Sciences \\
 52, B. Tulskaya, 115191 Moscow, Russia \\ [4pt]
$^2$ North-Eastern Federal University \\
58, Belinskogo, 677000 Yakutsk, Russia \\ [4pt]
e-mail: vabishchevich@gmail.com
\hfill Received: December 18, 2014 \\ [12pt]

\end{document}